\documentclass[12pt,a4paper,twoside]{article}
\usepackage{amsfonts}
\usepackage{amsthm,amsmath,amssymb}
\usepackage{color}
\def\trait #1 #2 #3 {\vrule width #1pt height #2pt depth #3pt}
\def\fin{
    \trait .3 5 0
    \trait 5 .3 0
    \kern-5pt
    \trait 5 5 -4.7
    \trait 0.3 5 0
\medskip}
\newtheorem{teor}{Theorem}[section]
\newtheorem{defin}[teor]{Definition}
\newtheorem{lemm}[teor]{Lemma}
\newtheorem{osse}[teor]{Remark}
\newtheorem{prop}[teor]{Proposition}
\newtheorem{defi}[teor]{Definition}
\newtheorem{coro}[teor]{Corollary}
\newtheorem{prob}[teor]{Problem}
\newtheorem{hypo}[teor]{Hypothesis}
\newcommand{\bele}{\begin{lemm}\begin{sl}}
\newcommand{\enle}{\end{sl}\end{lemm}}
\newcommand{\bedef}{\begin{defi}\begin{sl}}
\newcommand{\eddef}{\end{sl}\end{defi}}
\newcommand{\bete}{\begin{teor}\begin{sl}}
\newcommand{\ente}{\end{sl}\end{teor}}
\newcommand{\beos}{\begin{osse}\begin{rm}}
\newcommand{\eddos}{\end{rm}\end{osse}}
\newcommand{\bepr}{\begin{prop}\begin{sl}}
\newcommand{\empr}{\end{sl}\end{prop}}
\newcommand{\bepro}{\begin{prob}\begin{rm}}
\newcommand{\empro}{\end{rm}\end{prob}}
\newcommand{\bede}{\begin{defin}\begin{sl}}
\newcommand{\edde}{\end{sl}\end{defin}}
\newcommand{\beco}{\begin{coro}\begin{sl}}
\newcommand{\enco}{\end{sl}\end{coro}}
\newcommand{\behy}{\begin{hypo}\begin{sl}}
\newcommand{\enhy}{\end{sl}\end{hypo}}

\newcommand{\RR}{\mathbb{R}}

\newcommand{\beeq}[1]{\begin{equation}\label{#1}}
\newcommand{\eddeq}{\end{equation}}
\newcommand{\beeqa}[1]{\begin{eqnarray}\label{#1}}
\newcommand{\eddeqa}{\end{eqnarray}}
\newcommand{\beal}[1]{\begin{align}\label{#1}}
\newcommand{\eddal}{\end{align}}
\newcommand{\bespl}[1]{\begin{split}\label{#1}}
\newcommand{\edspl}{\end{split}}
\newcommand{\bega}[1]{\begin{gather}\label{#1}}
\newcommand{\edga}{\end{gather}}
\newcommand{\beeqax}{\begin{eqnarray*}}
\newcommand{\eddeqax}{\end{eqnarray*}}
\def\qed{\ifmmode   \else \leavevmode\unskip\penalty9999 \hbox{}\nobreak\hfill
  \fi
  \quad\hbox{\hskip.5em\vrule width.4em height.6em depth.05em\hskip.1em}}
\def\endproofsym{\qed}

\def\endnobox{\def\endproofsym{}\end{proof}\def\endproofsym{\qed}}

\newcommand{\no}{\nonumber}
\newcommand{\beeqao}{\begin{eqnarray}\no}
\newcommand{\bealo}{\begin{align}\no}
\newcommand{\besplo}{\begin{split}\no}
\newcommand{\begao}{\begin{gather}\no}
\def\trait #1 #2 #3 {\vrule width #1pt height #2pt depth #3pt}
\def\fin{\hfill
    \trait .3 5 0
    \trait 5 .3 0
    \kern-5pt
    \trait 5 5 -4.7
    \trait 0.3 5 0
\medskip}

\newcommand{\vc}[1]{{\boldsymbol #1}}
\newcommand{\dt}{\partial_t}

\newcommand{\io}{\int_\Omega}

\newcommand{\bd}{\boldsymbol{d}}

\newcommand{\rhs}{right hand side}

 \DeclareMathOperator{\dive}{div}

\let\TeXchi\chi
\def\chi{{\setbox0 \hbox{\mathsurround0pt
$\TeXchi$}\hbox{\raise\dp0 \copy0 }}}

\newcommand{\bu}{{\bf u}}

\newcommand{\Grad}{\nabla_x}
\newcommand{\tn}[1]{\mbox {\F #1}}

\font\F=msbm10   
scaled 1000

\newcommand{\db}{\boldsymbol{d}}

\newcommand{\ub}{{\bf u}}
\newcommand{\vb}{{\bf v}}

\newcommand{\teta}{\theta}

\def\fine{\hfill\kern4pt \vrule height4pt depth0pt width4pt }

\def\dive{\mbox{\rm div\,}}

scaled \magstep1 

   \numberwithin{equation}{section}
\begin{document}

\title{A new approach to non-isothermal models\\ for nematic liquid crystals}

\author{
Eduard Feireisl\thanks{Institute of Mathematics of the Czech
Academy of Sciences, \v Zitn\' a 25, 115 67 Praha 1, Czech Republic. E-mail:
{\tt feireisl@math.cas.cz} \ . The work of E.F. was supported by Grant
201/09/0917 of GA \v CR in the framework of research programmes
supported by AV\v CR Institutional Research Plan
AV0Z10190503}\\
\\
Michel Fr\'emond\thanks{Dipartimento di Ingegneria Civile, Universit\`a di Roma ``Tor Vergata'',
Via del Politecnico, 1, I-00133 Roma, Italy. E-mail:
{\tt Michel.Fremond@uniroma2.it}}\\
\\
Elisabetta Rocca\thanks{Dipartimento di Matematica, Universit\`a degli Studi di Milano,
Via Saldini 50, 20133 Milano, Italy. E-mail {\tt elisabetta.rocca@unimi.it}\ . The work of E.R. was partially supported by the FP7-IDEAS-ERC-StG
Grant \#256872 (EntroPhase) and by the MIUR-PRIN Grant 20089PWTPS
``Mathematical Analysis for inverse problems towards
applications''}\\
\\
Giulio Schimperna\thanks{Dipartimento di Matematica, Universit\`a degli Studi di Pavia,
Via Ferrata 1, 27100 Pavia,  Italy. E-mail: {\tt giusch04@unipv.it}\ . The work of G.S. was supported by the MIUR-PRIN Grant
2008ZKHAHN ``Phase transitions, hysteresis and multiscaling''}
\\
}
\date{}

\maketitle

\begin{abstract}
\noindent
We introduce a new class of
non-isothermal models describing the evolution of nematic liquid crystals
and prove their consistency with the fundamental laws of classical Thermodynamics.
The resulting system of equations captures all
essential features of  physically relevant models, in particular,
the effect of stretching of the director field is taken into account.
In addition, the associated initial-boundary value problem admits
global-in-time weak solutions without any essential
restrictions on the size of the initial data.
\end{abstract}

\section{Introduction}
\label{sec:intro}

The celebrated Leslie-Ericksen model of liquid crystals, introduced  by  Ericksen \cite{Er}
and Leslie \cite{Le78}, is a system of partial differential equations coupling
the Navier-Stokes equations governing the time evolution of the fluid velocity $\bu = \bu(t,x)$ with a Ginzburg-Landau
type equation describing the motion of the director field $\bd = \bd(t,x)$, representing
preferred orientation of molecules in a neighborhood of any point of a reference domain.

A considerably simplified version of the Leslie-Ericksen model was proposed by Lin and Liu
\cite{LinLiusimply}, \cite{LL01}, and subsequently analyzed by many authors, see \cite{SW10}, \cite{S02}, \cite{WXL10}
among others. The simplified model ignores completely the stretching and rotation effects of the director field induced by the straining of the fluid, which can be viewed
as a serious violation of the underlying physical principles.

Such a stretching term was subsequently treated by Coutand and Shkoller \cite{CoutShkoller}, who proved a local well-posedness
result for the corresponding model without thermal effects. The main peculiarity of this model is that the presence of
the stretching term causes the loss of the the total energy balance, which, indeed, ceases
to hold. In order to prevent this failure, Sun and Liu \cite{sunliu} introduced a variant of the
model proposed by Lin and Liu, where the stretching term is included in the system and a
new component added to the stress tensor in order to save the total energy balance.
A more general class of models based on the so-called Q-tensor formulation was recently
introduced in \cite{BZ,MZ} in the isothermal case.

Motivated by these considerations, in the the present
contribution, we propose a new approach to the modeling of
non-isothermal liquid crystals, based on the principles of
classical Thermodynamics and accounting for stretching and
rotation effects of the director field. To this end, we
incorporate the dependence on temperature into the model,
obtaining a complete energetically closed system, where the total
energy is conserved, while the entropy is being produced as the
system evolves in time. We apply here the mechanical methodology
of \cite{Frem}, which basically consists in deriving the equations
of the model by means of a generalized variational principle. This
states that the free energy $\Psi$ of the system, depending on the
proper {\sl state variables}, tends to decrease in a way that is
prescribed by the expression of a second functional, called
pseudopotential of dissipation, that depends (in a convex way) on
a set of {\sl dissipative variables}. In this approach, the stress
tensor $\sigma$, the density of energy vector ${\bf B}$ and the
energy flux tensor $\tn{H}$ are decoupled into their {\sl
non-dissipative}\/ and {\sl dissipative}\/ components, whose
precise form is prescribed by proper constitutive equations (see
below for details). It is interesting to note that the form of the
extra stress in the Navier-Stokes system obtained by this method
coincides with the formula derived from different principles by
Sun and Liu in \cite{sunliu}.

The system of partial differential equations resulting from this approach
couples the incompressible Navier-Stokes system
for the velocity $\ub$, with a Ginzburg-Landau type equation for the director field $\bd$ and a
total energy balance together with an entropy inequality, governing the dynamics of the absolute temperature $\theta$
of the system.

Leaving to the next Section~\ref{sec:model} the complete derivation of the model,
let us just briefly introduce here the PDE system we deal with.
The Navier-Stokes system couples  the incompressiblity condition
\begin{equation} \label{incompressiINTRO}
\dive \bu = 0
\end{equation}
with the conservation of momentum
\begin{equation} \label{momentumbalINTRO}
\bu_t + \bu \cdot \Grad \bu + \Grad p = \dive \tn{S} + \dive \sigma^{nd} + \vc{g},
\end{equation}
where $p$ is the pressure, and the stress is decomposed in a dissipative and non dissipative part,
respectively given by
\begin{equation}\label{stressINTRO}
\tn{S} =  \frac{\mu(\theta)}{2} \left( \Grad \bu + \Grad^t \bu \right),\
\sigma^{nd} = -\lambda \Grad \bd \odot\Grad\bd + \lambda(\vc{f}(\bd) - \Delta \bd )  \otimes \bd,
\end{equation}
where we have set $ \Grad \bd \odot \Grad\bd:=\sum_k \partial_i d_k \partial_j d_k$.

The director field equation has the form
\begin{equation} \label{micromouvINTRO}
\bd_t + \bu \cdot \Grad \bd - \bd \cdot \Grad \bu  = \gamma( \Delta \bd - \vc{f}(\bd)),
\end{equation}
where $\vc{f}(\bd)=\partial_{\bd} F(\bd)$ and $F$ penalizes the
deviation of the length $|\bd|$ from the value 1. It is a quite general function of $\bd$ that can be written
as a sum of a convex (possibly non smooth) part, and a smooth, but possibly non-convex one.
A typical example is $F(\bd)=(|\bd|^2-1)^2$.

Finally, the total energy balance
\begin{equation} \label{totaleINTRO}
\partial_t \left( \frac{1}{2}|\bu|^2 + e \right)
+ \bu \cdot \Grad \left( \frac{1}{2}|\bu|^2 + e \right)
+ \dive \Big( p \bu + \vc{q} - \tn{S} \bu - \sigma^{nd} \bu \Big)
\end{equation}
\[
= \vc{g} \cdot \bu + \lambda \gamma \dive \Big( \Grad \bd \cdot \left( \Delta \bd - \vc{f}(\bd) \right) \Big),
\]
with the internal energy and the flux
\[
e = \frac{\lambda}{2} |\Grad \bd |^2 + \lambda F(\bd) + \theta,
\quad {\bf q}= {\bf q}^d- \lambda \Grad \bd \cdot \Grad \bu \cdot \bd, \quad {\bf q}^d=-k(\theta )\Grad\theta -h(\theta )(\db\cdot \Grad\theta )\db,
\]
is coupled with the entropy inequality
\begin{equation} \label{eqtetaINTRO}
H(\theta)_t + \bu \cdot \Grad H(\theta) + \dive (H'(\theta) \vc{q}^d )
\end{equation}
\[
\geq
H'(\theta) \left( \tn{S} : \Grad \bu + \lambda \gamma | \Delta \bd - \vc{f}(\bd) |^2 \right) +
H''(\theta) \vc{q}^d \cdot \Grad \theta,
\]
holding true for any smooth non-decreasing concave function $H$. The derivation of
the above system will be detailed in the next Section~\ref{sec:model}, while the
remainder of the paper will be devoted to the proof of a global existence
result for the corresponding initial-boundary value problem in the framework of
weak solutions in $\Omega\times (0,T)$, being
$\Omega$ a bounded and sufficiently regular subset of $\RR^3$ and $T$ a given final time.

Let us note that the model obtained here looks quite different from the one obtained in
\cite{frs}. This is mainly due to the presence
in the internal energy $e$ of the quadratic term $|\Grad\db|^2$ (that is related to
the expression \eqref{psi} of the free energy functional) and
to the stretching term $\bd\cdot\Grad\ub$ in \eqref{micromouvINTRO} which produces, in order
that the principles of Thermodynamics are respected, two new non dissipative contributions in the
stress tensor $\tn{S}$ in \eqref{stressINTRO} and in the flux ${\bf q}$.
Actually, the latter is given here by the sum of a standard heat flux and
of an elastic part given by the term $- \lambda \Grad \bd \cdot \Grad \bu \cdot \bd$
(cf.\ the next Section~\ref{sec:model} for further details on this point).

Indeed, in contrast with \cite{frs}, the  presence of the stretching
term $\bd\cdot\Grad\ub$ in the director field equation prevent us from applying
any form of the maximum principle to \eqref{micromouvINTRO}.
Hence, we cannot recover an $L^\infty$-bound on $\bd$ (which we obtained, instead, in \cite{frs}).
However, we can still get here the global existence of weak solutions to the initial boundary
value problem coming from the PDE system (\ref{incompressiINTRO}--\ref{eqtetaINTRO}) without imposing any restriction
on the space dimension, on the size of the initial data or on the viscosity coefficient $\mu$
(such a restriction was taken in the paper \cite{sunliu}, devoted to an {\sl isothermal}\/
model closely related to ours). In this sense, our results can be seen as a generalization of
those obtained in \cite{sunliu}.

The compatibility of the model
with \emph{First and Second laws} of thermodynamics turns out to be the main source of {\it a priori} bounds
that can be used, in combination with compactness arguments, to ensure stability of the family of approximate solutions.
The key point of this approach is replacing the heat equation, commonly used in models of heat conducting fluids,
by the total energy balance \eqref{totaleINTRO}. Accordingly, the resulting
system of equations is free of dissipative terms that are difficult to handle,
due to the low regularity of the weak solutions. In contrast with the standard
theory of Navier-Stokes equations, however, we have to control the
pressure appearing explicitly in the total energy flux and, in order to do that we will
need to assume the complete slip boundary conditions on the velocity filed $\ub$ (cf. \eqref{slip}).
Note that a similar method applied to different models has been
recently used in \cite{BFM}, \cite{FM06}, and \cite{frs}.

Finally, let us notice that the non-isothermal liquid crystal model accounting for the stretching contribution
has also recently been analyzed in \cite{cr} (in case of Dirichlet boundary conditions for $\ub$ and Neumann or
non homogeneous Dirichlet boundary conditions for $\bd$) and in \cite{prs}, where the long time behaviour of solutions
is investigated in two cases: in the 3D case without any condition on the size of the viscosity coefficient $\mu$
and in case of a non analytic nonlinearity $\vc{f}$. Both these results generalize the ones obtained
in  \cite{WXL10}.

The paper is organized as follows.
In Section~\ref{sec:model}, we give a detailed derivation
of the model and discuss its compatibility with the basic laws of
Thermodynamics. In Section~\ref{sec:mainres}, we introduce some technical hypotheses and formulate
the main result concerning existence of global-in-time weak solutions to the resulting PDE system.
Finally, the last two Sections~\ref{a} and \ref{A} are
devoted to the proof of the existence result via approximation, a-priori estimates and passage to the limit techniques
based on lower semicontinuity and convexity arguments. As already pointed out, the energy balance is written in the
form of a conservation law for the total energy rather than for the temperature, where the highly non-linear terms
dissipative terms are absent. The price to pay is the explicit appearance of the \emph{pressure} in the global energy
balance determined implicitely by the Navier-Stokes system.


\section{Mathematical model}
\label{sec:model}

We suppose that the fluid occupies a bounded spatial domain $\Omega
\subset \RR^3$, with a sufficiently regular boundary, and denote by
$\bu = \bu(t,x)$ the associated {\em velocity field} in the Eulerian reference system.
Moreover, we introduce the {\em absolute temperature}  $\theta(t,x)$ and the {\em director field}
$\db (t,x)$, representing the
preferred orientation of molecules in a neighborhood of any point
of the reference domain. Furthermore,
we denote
$$
\frac{d w}{dt} = \dot{w} =w_t+\ub\cdot\Grad w,
$$
the {\em material derivative} of a generic function $w$, while
$w_t$ (or also $\dt w$) denotes the partial derivative with respect to $t$.

Finally, the quantity
\[
\frac {D \db }{Dt}  = \db_t + \ub \cdot \Grad \db - \db \cdot
\Grad \ub
\]
characterizes the total transport of the orientation vector $\db$. Note that the last term accounts for
stretching of the director field induced by the straining of the fluid.

\subsection{Free-energy and pseudopotential of dissipation}

Following the general approach proposed in the monograph \cite{Frem},
we start by specifying, in agreement with
the principles of classical Thermodynamics, the free-energy and the
pseudopotential of dissipation. The interested reader may consult \cite[Chapters 2,3]{Frem} for details.

We begin by introducing the set of {\em state variables}, describing the
actual configuration of the material, specifically,
\begin{equation*}
E=(\bd,\Grad \bd,\theta).
\end{equation*}
Next, the set of the {\em dissipative variables} describing the evolution of the system,
and, in particular, the way it dissipates energy,
is given by
\begin{equation*}
\delta E=\left(\varepsilon (\ub),\frac{ D \bd}{Dt}, \Grad\theta\right),
\end{equation*}
where
$$\varepsilon (\ub):=\frac{\left(\Grad \ub+\Grad^t \ub\right)}{2} $$
denotes the symmetric gradient of $\ub$.

Motivated by the original (isothermal) theory proposed by Ericksen \cite{Er61} and Leslie \cite{Le63},
we choose the free energy functional in the form
\begin{equation}\label{psi}
\Psi (E)=\frac{\lambda}{2}|\Grad\bd| ^{2}+ \lambda F(\bd)-\theta \log \theta\,,
\end{equation}
where $\lambda$ is a positive constant.
The function $F$ in \eqref{psi} penalizes the
deviation of the length $|\bd|$ from its natural value 1; generally, $F$ is assumed to be
a sum of a dominating convex (and possibly non smooth) part
and a smooth non-convex perturbation of controlled growth.
A typical example is $F(\bd)=(|\bd|^2-1)^2$. For the sake of simplicity, we have assumed
that the thermal and elastic effects are uncoupled in $\Psi$.

The evolution of the system is characterized by a second functional $\Phi$,
called {\em pseudopotential of dissipation}, assumed to be nonnegative
and convex with respect to the dissipative variables. Specifically, we consider
$\Phi$ in the form
\begin{align*}
\Phi (\delta E,E)=&\,\frac{\mu (\theta )}{2}|\varepsilon (\ub)|^{2}+I_{0}(\dive \ub)+\frac{
k(\theta )}{2\theta }|\Grad\theta |^{2}+\frac{\eta}{2}\left|\frac{D\bd}{Dt}\right|^{2}+\frac{
h(\theta )}{2\theta }|\bd\cdot \Grad\theta |^{2}\,,
\end{align*}
where $\mu = \mu(\theta) > 0$ is the viscosity coefficient, $\eta > 0$ is a constant, and $k$, $h$ represent the heat conductivity
coefficients - positive functions of the temperature. The {\em incompressibility} of the fluid is formally
enforced by $I_0$ - the indicator function of $\{0\}$ (given
by $I_0 = 0$ if $\dive \ub =0$ and $+\infty$ otherwise).

\subsection{Constitutive relations}

We start by introducing the stress tensor
$\sigma$, the density of energy vector ${\bf B}$, and the energy flux tensor $\tn{H}$; all assumed to be the sum of their
non-dissipative and dissipative components, namely, $\sigma=\sigma^{nd}+\sigma^d$, ${\bf B}={\bf B}^{nd}+{\bf B}^d$, $\tn{H}=\tn{H}^{nd}+\tn{H}^d$,
where
\begin{gather}\label{Bnd}
{\bf B}^{nd}=\frac{\partial \Psi }{\partial \db}=  \lambda \frac{\partial
F}{\partial \db}=: \lambda \vc{f}(\bd),
\\
\label{Bd}
{\bf B}^{d}=\frac{\partial \Phi }{\partial \frac{D\db}{Dt}} = \eta \frac{D\db}{Dt},\\
\label{Hnd}
\tn{H}^{nd}=\frac{\partial \Psi }{\partial \Grad \db}=\lambda\Grad \db.
\end{gather}
Moreover, we set  $\tn{H}^d \equiv 0$.

The heat and entropy fluxes (denoted respectively by ${\bf q}^d$ and ${\bf Q}$) are
\begin{equation}\label{Q}
{\bf q}^d=\theta {\bf Q} = -\theta \frac{\partial \Phi }{\partial \Grad\theta }%
=-k(\theta )\Grad\theta -h(\theta )(\db\cdot \Grad\theta )\db.
\end{equation}

The stress tensor $\sigma$ consists of two parts: the dissipative one
\begin{equation}\label{sigmad}
\sigma ^{d}=\frac{\partial \Phi }{\partial \varepsilon (\ub)}=\mu (\theta
)\varepsilon (\ub)-p\mathbb{I}=:\tn{S}-p\tn{I}\,,
\end{equation}
\[
- p \in \partial I_0 (\dive \bu), \ \tn{S} = \mu(\theta) \varepsilon (\ub),
\]
and the non dissipative part $\sigma^{nd}$ to be determined below
(cf.~\eqref{eqs} and \eqref{eqe}).

The entropy of the system is given by
\begin{equation}\label{s}
s=-\frac{\partial \Psi }{\partial \theta }=1+\log \theta
\end{equation}
and, finally, the internal energy $e$ reads
\begin{equation}\label{e}
e=\Psi +\theta s=\theta +\lambda F(\bd)
+ \frac{\lambda}{2} {|\Grad\bd|^2}.
\end{equation}

\subsection{Field equations}

In accordance with \emph{Newton's second law}, the balance of momentum reads
\begin{equation}\label{mombal}
\partial_t \ub  + \dive (\ub \otimes \ub )  = \dive \sigma + \vc{g} \, ,
\end{equation}
where $\vc{g}$ is a given external force.

The entropy balance can be written in the form
\begin{equation}\label{eqs}
s_t+\ub\cdot\Grad s+\dive {\bf Q}=\frac{1}{\theta }\left( \sigma
^{d}:\varepsilon(\ub) +{\bf B}^{d}\cdot\frac{D \bd}{Dt}-{\bf
Q}\cdot\Grad\theta \right)
\end{equation}
or equivalently
\begin{equation}\label{eqsbis}
\teta \frac{d s}{dt}+\dive {\bf q}^d=\sigma
^{d}:\varepsilon(\ub) +{\bf B}^{d}\cdot\frac{D \bd}{D t}.
\end{equation}
In agreement with {\em Second law of Thermodynamics},
the right hand side of \eqref{eqs} is non-negative.

The balance of internal energy reads
\begin{equation}\label{eqe}
e_t+\ub\cdot\Grad e+\dive {\bf q}=\sigma :\varepsilon(\ub)+{\bf
B}\cdot\frac{D\bd}{Dt}+\tn{H}:\Grad\frac{D\bd}{Dt}\,,
\end{equation}
with the internal energy flux $\vc{q} = \vc{q}^d + \vc{q}^{nd}$, where the dissipative
part $\vc{q}^d$ is given by (\ref{Q}), while the non-dissipative component will be determined below.

Finally, the equation which rules the evolution of the orientation vector $\bd$ is derived from the
principle of virtual powers (cf. \cite[Chap. 2]{Frem}) and it takes the form
\begin{equation} \label{pom}
\dive \tn{H} -{\bf B}={\bf 0}\,,
\end{equation}
specifically,
\begin{equation} \label{eqd}
  \db_t + \ub \cdot \Grad \db - \db \cdot \Grad \ub
   = {\gamma} (\Delta \bd - \vc{f}(\bd)), \ \gamma = \lambda / \eta.
\end{equation}

The non-dissipative component of the stress $\sigma^{nd}$ and of the flux
$\vc{q}^{nd}$ are determined by means of \eqref{eqs}, \eqref{eqe}, and the constitutive relations derived above.
Indeed, computing $\frac{d e}{dt}$ by means of the standard Helmholtz relations, we get
\begin{equation}\label{ide}
\frac{ d e}{dt}= \frac{d\Psi}{dt}+\theta\frac{d s}{dt}+\frac{d \theta}{dt} s=\Psi_{\bd}\cdot\frac{d\bd}{dt}
+\Psi_{\Grad\bd}:\frac{d(\Grad\bd)}{dt}+\theta\frac{d s}{dt} ,
\end{equation}
whereas
\begin{equation}\label{psigradb}
\Psi_{\Grad\bd} : \frac{d(\Grad\bd)}{dt}=\tn{H}^{nd}:\left( \Grad\frac{d\bd}{dt}-\Grad\ub \cdot \Grad\bd\right).
\end{equation}
Thus, rewriting \eqref{ide} with help of \eqref{psigradb}, and expressing $\theta\frac{d s}{dt}$ by means of \eqref{eqsbis},
we get, thanks also to (\ref{Bnd}--\ref{sigmad}),
\begin{equation}\label{edot}
\dive \vc{q}^{nd} - \sigma^{nd} : \Grad \bu = - \lambda (\vc{f}(\bd) \otimes \bd) : \Grad \bu - \lambda \sum_{i,j,k}
\partial_{x_k} d_i (\partial^2_{x_j,x_k} u_i ) d_j .
\end{equation}
Therefore,
\begin{equation}\label{sigmand}
\sigma^{nd}=-\lambda \Grad \bd\odot\Grad\bd + \lambda(\vc{f}(\bd) - \Delta \bd )  \otimes \bd,\
\vc{q}^{nd} = - \lambda \Grad \bd \cdot \Grad \bu \cdot \bd .
\end{equation}

Summing up the previous discussion,
we arrive at the following system of equations:

\centerline{\textsc{incompressibility:}}
\begin{equation} \label{incompressi}
\dive \bu = 0;
\end{equation}
\centerline{\textsc{conservation of momentum:}}
\begin{equation} \label{momentumbal}
\bu_t + \bu \cdot \Grad \bu + \Grad p = \dive \tn{S} + \dive \sigma^{nd} + \vc{g},
\end{equation}
where $p$ is the pressure, and
\begin{equation}\label{stress}
\tn{S} =  \frac{\mu(\theta)}{2} \left( \Grad \bu + \Grad^t \bu \right),\
\sigma^{nd} = -\lambda \Grad \bd \odot\Grad\bd + \lambda(\vc{f}(\bd) - \Delta \bd )  \otimes \bd;
\end{equation}
\centerline{\textsc{director field equation:}}
\begin{equation} \label{micromouv}
\bd_t + \bu \cdot \Grad \bd - \bd \cdot \Grad \bu  = \gamma( \Delta \bd - \vc{f}(\bd)),
\end{equation}
\centerline{\textsc{total energy balance:}}
\begin{equation} \label{totale}
\partial_t \left( \frac{1}{2}|\bu|^2 + e \right)
+ \bu \cdot \Grad \left( \frac{1}{2}|\bu|^2 + e \right)
+ \dive \Big( p \bu + \vc{q} - \tn{S} \bu - \sigma^{nd} \bu \Big)
\end{equation}
\[
  = \vc{g} \cdot \bu + \lambda \gamma \dive \Big( \Grad \bd \cdot \left( \Delta \bd - \vc{f}(\bd) \right) \Big),
\]
with the internal energy
\[
  e = \frac{\lambda}{2} |\Grad \bd |^2 + \lambda F(\bd) + \theta
\]
and the flux
\[
  {\bf q} = {\bf q}^d+{\bf q}^{nd}
  = -k( \theta )\Grad\theta -h(\theta )(\db\cdot \Grad\theta )\db
  - \lambda \Grad \bd \cdot \Grad \bu \cdot \bd,
\]
together with

\centerline{\textsc{entropy inequality:}}
\begin{equation} \label{eqteta}
  H(\theta)_t + \bu \cdot \Grad H(\theta) + \dive (H'(\theta) \vc{q}^d )
\end{equation}
\[
\geq
H'(\theta) \left( \tn{S} : \Grad \bu + \lambda \gamma | \Delta \bd - \vc{f}(\bd) |^2 \right) +
H''(\theta) \vc{q}^d \cdot \Grad \theta,
\]
holding for any smooth, non-decreasing and concave function $H$.

Actually, the total energy balance (\ref{totale}) follows easily from
(\ref{mombal}), (\ref{eqe}), combined with (\ref{pom}).
It is remarkable that equations (\ref{incompressi}--\ref{micromouv})
in the isothermal case reduce to the model derived by Sun and Liu in
\cite{sunliu} by means of a different method.


\section{Main results}
\label{sec:mainres}


\subsection{Initial and boundary conditions}

In view of a rigorous mathematical study, system (\ref{incompressi}--\ref{micromouv})
must be supplemented by suitable boundary conditions.
Actually, to avoid the effect of boundary layer on the motion,
we assume \emph{complete slip} boundary conditions
for the velocity:
\begin{equation}\label{slip}
\ub \cdot \vc{n}|_{\partial \Omega} = 0,\ [(\tn{S} + \sigma^{nd}) \vc{n}] \times \vc{n} |_{\partial \Omega} = 0.
\end{equation}
Moreover, we consider \emph{no-flux} boundary condition for the temperature
\begin{equation} \label{noflux}
\vc{q}^d \cdot \vc{n}|_{\partial \Omega} = 0,
\end{equation}
and Neumann boundary condition for the director field
\begin{equation}\label{neumann}
\Grad d_i \cdot \vc{n}|_{\partial \Omega} = 0 \ \mbox{for}\ i=1,2,3.
\end{equation}
The last relation accounts for the fact that there is no contribution to the surface
force from the director $\bd$. Note that the above conditions are also suitable
for the implementation of a numerical scheme (see \cite{LS} for further comments
on this point).

Of course, we also need to assume the initial conditions
\begin{equation} \label{initial}
\vc{u}(0, \cdot) = \vc{u}_0, \ \vc{d}(0, \cdot) = \vc{d}_0, \
\teta(0,\cdot) = \teta_0,
\end{equation}
In the remaining part of the paper, our aim will be that of showing
existence of global-in-time solutions to system
(\ref{incompressi}--\ref{eqteta}),
coupled with the above initial and boundary conditions and
without assuming any essential restriction on the data.


\subsection{Weak formulation}
\label{weak}

In the weak formulation, the momentum equation (\ref{momentumbal}), together with
the incompressibility constraint (\ref{incompressi}),
and the boundary conditions (\ref{slip}), are replaced by a family of integral identities
\begin{equation}\label{weak1}
\int_{\Omega} \vc{u}(t, \cdot) \cdot \Grad \varphi = 0\ \ \mbox{for a.a.} \ t \in (0,T)
\end{equation}
for any test function $\varphi \in C^{\infty} (\overline{\Omega})$,
and
\begin{equation} \label{weak2}
\int_0^T \int_{\Omega} \Big( \vc{u} \cdot \partial_t \varphi + \vc{u} \otimes \vc{u} : \Grad \varphi +
p \, \dive \varphi \Big)
\end{equation}
\[
=
\int_0^T \int_{\Omega} (\tn{S} + \sigma^{nd}) : \Grad \varphi - \int_\Omega
\vc{g} \cdot \varphi - \int_{\Omega} \vc{u}_0 \cdot \varphi (0, \cdot)\,,
\]
for any $\varphi \in C^\infty_0 ([0,T) \times \overline{\Omega}; \RR^3)$,
$\varphi \cdot \vc{n}|_{\partial \Omega} = 0$. Note that (\ref{weak2}) includes also the initial
condition $\bu(0, \cdot) = \bu_0$.

Equation (\ref{micromouv}) describing the evolution of the director field $\bd$ will be satisfied in the strong sense,
more specifically,
\begin{equation} \label{weak3}
\partial_t \bd + \ub\cdot\Grad \bd-\bd\cdot\Grad \ub = \gamma \Big( \Delta \bd - \vc{f}(\bd) \Big)\ \mbox{ a.e. in }(0,T) \times \Omega,
\end{equation}
together with
\[
\quad \Grad\bd_i\cdot \vc{n}_{|\partial\Omega}=0\,, \ i=1,2,3,\quad \bd (0,\cdot)=\bd_0.
\]

Similarly, the weak formulation of the total energy balance (\ref{totale}) reads
\begin{equation} \label{weak4}
\int_0^T \int_{\Omega} \left(\left( \frac{1}{2} |\vc{u}|^2
+ e \right) \partial_t \varphi\right)+\int_0^T \int_{\Omega} \left(\left(\frac{1}{2} |\vc{u}|^2 + e \right) \vc{u} \cdot \Grad
\varphi \right)
\end{equation}
\[
+ \int_0^T \int_{\Omega} \left(p \bu + \vc{q} - \tn{S} \bu - \sigma^{nd} \bu \right)\cdot \Grad \varphi
\]
\[=
\lambda\gamma\int_0^T \int_{\Omega}\Big( \Grad \bd \cdot \left( \Delta \bd - \vc{f}(\bd) \right) \Big)
\cdot \Grad \varphi -\int_0^T \int_{\Omega} \vc{g}\cdot \ub \varphi  - \int_{\Omega}
\left( \frac{1}{2} |\vc{u}_0|^2 + e_0 \right)
\varphi (0, \cdot)\,,
\]
for any $\varphi \in C^{\infty}_0 ([0,T) \times \overline{\Omega})$,
where $e_0=\frac{\lambda}{2}|\Grad\bd_0|^2+\lambda F(\bd_0)+\teta_0$.

Finally, the entropy inequality \eqref{eqteta} is replaced by
\begin{equation}\label{weak5}
\int_0^T\io H(\teta) \partial_t\varphi +\int_0^T\io \left(H(\teta) \ub +H'(\teta)\vc{q}^d\right)\cdot\Grad\varphi
\end{equation}
\[
\leq -\int_0^T\io \left(H'(\teta)\left(\tn{S}: \Grad \bu + \lambda \gamma | \Delta \bd - \vc{f}(\bd) |^2 \right) +
H''(\theta) \vc{q}^d \cdot \Grad \theta\right)\varphi-\io H(\teta_0)\varphi(0,\cdot)
\]
for any $\varphi \in C^{\infty}_0 ([0,T) \times \overline{\Omega})$, $\varphi\geq 0$, and for any smooth, non-decreasing and concave function
$H$.

A {\em weak solution} is a triple $(\ub,\,\bd,\,\teta)$ satisfying
(\ref{weak1}--\ref{weak5}).

\subsection{Main existence theorem}

Before formulating the main result of this paper, we
list the hypotheses imposed on the constitutive functions. Specifically, we assume
that
\begin{equation}\label{hyp1}
F \in C^2(\RR^3), \quad F \geq 0,\quad F \ \mbox{ convex for all}\
 |\vc{d}| \geq D_0, \ \lim_{|\bd|\to\infty} F(\bd)=\infty,
\end{equation}
for a certain $D_0>0$.

The transport coefficients $\mu$, $k$, and $h$ are continuously
differentiable functions of the absolute temperature satisfying
\begin{equation} \label{hyp2}
0 < \underline{\mu} \leq \mu(\teta) \leq \overline{\mu}, \quad 0 <
\underline{k} \leq k(\teta),\,\, h
(\teta) \leq \overline{k} \ \mbox{ for all}\ \teta \geq 0
\end{equation}
for suitable constants $\underline{k}$, $\overline{k}$,
$\underline{\mu}$, $\overline{\mu}$.

Our main result reads as follows.

\bete \label{theo1} Let $\Omega \subset \RR^3$ be a bounded domain
of class $C^{2 + \nu}$ for some $\nu>0$, $\vc{g}\in L^2((0,T)\times\Omega;\RR^3)$. Assume that hypotheses (\ref{hyp1}),
(\ref{hyp2}) are satisfied. Finally, let the initial data be such
that
\begin{equation}\label{hyp4}
\begin{array}{c}
\vc{u}_0 \in L^2(\Omega; \RR^3),\ \dive \vc{u}_0 = 0,\ \vc{d}_0 \in
W^{1,2} (\Omega; \RR^3), \ F(\bd_0)\in L^1(\Omega),\\\\
 \teta_0 \in
L^1(\Omega),\ {\rm ess} \inf_{\Omega} \teta_0 > 0 .
\end{array}
\end{equation}

Then, problem (\ref{weak1}--\ref{weak5}) possesses a weak
solution $(\vc{u},\vc{d},\teta)$ in $(0,T) \times \Omega$
belonging to the class
\begin{equation} \label{reg1}
\vc{u} \in L^\infty(0,T; L^2(\Omega; \RR^3)) \cap L^2(0,T;
W^{1,2}(\Omega; \RR^3)),
\end{equation}
\begin{equation} \label{reg2}
\vc{d} \in  L^\infty(0,T;
W^{1,2}(\Omega; \RR^3))\cap L^2(0,T; W^{2,2}(\Omega;\RR^3)),
\end{equation}
\begin{equation}
\label{regoF}
F(\bd)\in L^\infty(0,T; L^1(\Omega))\cap L^{5/3}((0,T) \times \Omega),
\end{equation}
\begin{equation} \label{reg3}
\teta \in L^\infty(0,T; L^1(\Omega)) \cap L^p(0,T;
W^{1,p}(\Omega)),\ 1 \leq p < 5/4, \ \teta > 0 \ \mbox{a.e. in}\
(0,T) \times \Omega,
\end{equation}
with the pressure $p$,
\begin{equation} \label{reg4}
p \in L^{5/3}((0,T) \times \Omega).
\end{equation}
\ente

The rest of the paper is devoted to the proof of Theorem
\ref{theo1}.

\section{A priori bounds}
\label{a}

In this section, we collect the
available {\em a priori}\ estimates. These will assume a
rigorous character in the framework of the approximation scheme presented in
Section \ref{A} below.

Integrating (\ref{totale}) over $\Omega$ and using Gronwall's lemma, we immediately obtain
the following bounds:
\begin{equation}\label{apr1}
\vc{u} \in L^\infty(0,T; L^2(\Omega; \RR^3)),
\end{equation}
\begin{equation}\label{apr2}
\teta \in L^\infty(0,T; L^1(\Omega)),
\end{equation}
\begin{equation} \label{apr3}
\vc{d} \in L^\infty (0,T; W^{1,2}(\Omega; \RR^3)), \ F(\bd)\in L^\infty(0,T; L^1(\Omega)),
\end{equation}
where we have used hypotheses (\ref{hyp1}), (\ref{hyp2}).

Similarly, integrating \eqref{eqteta} with $H(\teta)=\teta$, and using \eqref{apr2}, we obtain
\begin{equation}\label{esDeltaf}
\varepsilon(\ub)\in L^2((0,T)\times \Omega, \RR^{3 \times 3}), \ \Delta \bd-\vc{f}(\bd)\in L^2((0,T)\times \Omega;
\RR^3).
\end{equation}
yielding, by virtue of \eqref{apr1} and Korn's inequality,
\begin{equation}\label{apr1bis}
 \vc{u} \in L^2(0,T;W^{1,2}(\Omega; \RR^3))\cap L^{10/3}((0,T)\times\Omega;\RR^3).
\end{equation}

Moreover, it follows from \eqref{esDeltaf} and convexity of $F$ (cf. hypothesis \eqref{hyp1}) that
\begin{equation}\label{apr4}
\vc{f}(\bd)\in L^2((0,T)\times \Omega; \RR^3);
\end{equation}
therefore, using \eqref{esDeltaf} again we infer that
\begin{equation} \label{apr5}
\vc{d} \in L^2(0,T; W^{2,2}(\Omega; \RR^3)).
\end{equation}

Interpolating \eqref{apr3} and \eqref{apr5} we get
\begin{equation}\no
\bd\in L^{10} ((0,T)\times \Omega; \RR^3), \ \Grad\bd\in L^{10/3}((0,T)\times \Omega; \RR^{3\times 3}),
\end{equation}
whence (cf.~\eqref{stress})
\begin{equation}\label{regsigmand}
\sigma^{nd}\in L^{5/3}((0,T)\times \Omega;\RR^{3 \times 3}).
\end{equation}
By the same token, by means of convexity of $F$ (cf. \eqref{hyp1}), we have
\[
  |F(\bd)|\leq c(1+|\vc{f}(\bd)||\bd|),
\]
yielding
\begin{equation}\label{regF}
  F(\bd)\in L^{5/3}((0,T)\times \Omega).
\end{equation}

As the velocity satisfies the slip boundary conditions (\ref{slip}),
the pressure $p$ can be ``computed'' directly from
(\ref{momentumbal}) as the unique solution of the elliptic problem
\begin{equation}\nonumber
\Delta p = \dive \dive \Big( \tn{S} +\sigma^{nd}- \vc{u} \otimes \vc{u} \Big)+\dive \vc{g},
\end{equation}
supplemented with the boundary condition
\begin{equation}\nonumber
\Grad p\cdot {\bf n} =\left(\dive{\left(\tn{S} +\sigma^{nd}- \vc{u} \otimes \vc{u}\right)+\vc{g} }\right)\cdot{\bf n}  \mbox{ on }\partial\Omega\,.
\end{equation}
To be more precise, the last two relations have to be interpreted in a ``very weak'' sense.
Namely, the pressure $p$ is determined through a family of integral identities:
\begin{equation} \label{press}
\int_\Omega p \Delta \varphi = \int_{\Omega} \Big(\tn{S} +\sigma^{nd} - \vc{u} \otimes
\vc{u} \Big) : \Grad^2 \varphi-\int_\Omega \vc{g}\cdot \Grad \varphi,
\end{equation}
for any test function $\varphi \in C^\infty(\overline{\Omega})$,
$\Grad \varphi \cdot \vc{n}|_{\partial \Omega} = 0$. Consequently,
the bounds established in (\ref{apr1bis}) and \eqref{regsigmand} may be used,
together with the standard elliptic regularity results, to conclude
that
\begin{equation}\label{apr7}
  p \in L^{5/3}((0,T) \times \Omega).
\end{equation}

Finally, the choice $H(\teta) = ( 1 + \teta )^\eta , \ \eta \in (0, 1)$, in
(\ref{eqteta}), together with the uniform bounds obtained in
(\ref{apr1}--\ref{apr1bis}), yields
\begin{equation} \label{apr8}
\Grad (1 + \teta)^{\nu} \in L^2((0,T) \times \Omega; \RR^3) \
\mbox{for any}\ 0 <\nu < \frac{1}{2}\,.
\end{equation}
Now, we apply an interpolation argument already exploited
in \cite{BFM}. Using (\ref{apr2}) and \eqref{apr8} and interpolating between
$\teta\in L^\infty(0,T;L^1(\Omega))$ and $\teta^\nu\in L^1(0,T;L^3(\Omega))$,
for $\nu\in (0,1)$,
we immediately  get
\begin{equation}\label{apr8bis}
 \teta\in {L^q((0,T) \times \Omega)} \ \mbox{for any}\ 1
  \leq q < 5/3\,.
\end{equation}
Furthermore, seeing that
\begin{equation}\nonumber
 \int_{(0,T)\times\Omega}|\Grad\teta|^p\leq
  \left(\int_{(0,T)\times\Omega}|\Grad\teta|^2\teta^{\nu-1}\right)^{\frac{p}{2}}
  \left(\int_{(0,T)\times\Omega}\teta^{(1-\nu)\frac{p}{2-p}}\right)^{\frac{2-p}{2}}\,
\end{equation}
for all $p\in [1,5/4)$ and $\nu>0$,
we conclude from \eqref{apr8} and \eqref{apr8bis} that
\begin{equation} \label{apr9}
  \Grad \teta \in L^p((0,T) \times \Omega; \RR^3) \ \mbox{for any}\ 1
   \leq p < 5/4.
\end{equation}

Finally, the same argument and $H(\teta) = \log\teta$ in
(\ref{eqteta}) give rise to
\begin{equation}\label{apr9bis}
 \log\teta \in L^2((0,T);W^{1,2}( \Omega))\cap L^\infty(0,T; L^1(\Omega)),
\end{equation}
where we have used \eqref{apr2}.

The {\it a priori}\ estimates derived in this section comply with
the regularity class (\ref{reg1}--\ref{reg4}). Moreover, it can be shown
that the solution set of (\ref{weak1}--\ref{weak5})
is weakly stable (compact) with respect to these
bounds, namely, any sequence of (weak) solutions that satisfies
the uniform bounds established above
has a subsequence that converges to some limit
that still solves the system.
Leaving the proof of weak sequential stability to the
interested reader, we pass directly to the proof of Theorem
\ref{theo1} constructing a suitable family of
\emph{approximate}\ problems.

\section{Approximations}

\label{A}

For the sake of simplicity, we restrict ourself to the case $\vc{g}=\vc{0}$ and $\lambda=\gamma=1$.
Solutions to the Navier-Stokes system (\ref{weak1}), (\ref{weak2})
will be constructed by means of the nowadays standard Faedo-Galerkin
approximation scheme, see Temam \cite{temam}. Let
$W^{1,2}_{n,\sigma}(\Omega; \RR^3)$ be the Sobolev space of
solenoidal functions satisfying the impermeability boundary
condition, specifically,
\begin{equation}\nonumber
W^{1,2}_{n,\sigma} = \{ \vc{v} \in W^{1,2}(\Omega; \RR^3) \ | \
\dive \vc{v} = 0 \ \mbox{a.e. in}\ \Omega,\ \vc{v} \cdot
\vc{n}|_{\partial \Omega} = 0 \}\,.
\end{equation}
Since $\partial \Omega$ is of class $C^{2+ \nu}$, there exists an
orthonormal basis $\{ \vc{v}_n \}_{n=1}^\infty$ of the Hilbert space
$W^{1,2}_{n , \sigma}$ such that $\vc{v}_n \in C^{2 + \nu}$, see
\cite[Theorem 10.13]{FN}.
We take $M\leq N$ and denote $X_N = {\rm span}\{ \vc{v}_n
\}_{n=1}^N$, and $[ \vc{v} ]_M$ -
the orthogonal projection onto the space ${\rm span}\{ \vc{v}_n
\}_{n=1}^M$.

The approximate velocity fields $\vc{u}_{N,M} \in C^1([0,T]; X_N)$ solve
the Faedo-Galer\-kin system
\begin{equation} \label{approx1}
 \frac{{\rm d}}{{\rm d}t} \int_\Omega \vc{u}_{N,M} \cdot \vc{v} =
  \int_{\Omega}  [ \vc{u}_{N,M} ]_M\otimes \vc{u}_{N,M} : \Grad \vc{v}
  - \frac{1}{M} \io |\Grad\vc{u}_{N,M}|^{r-2}\Grad\vc{u}_{N,M}:\Grad\vc{v}
\end{equation}
\begin{equation}\nonumber
 - \int_\Omega \frac{\mu(\teta_{N,M})}2
  \Big( \Grad \vc{u}_{N,M} + \Grad^t \vc{u}_{N,M} \Big): \Grad \vc{v}
 + \int_\Omega \Grad \vc{d}_{N,M} \odot \Grad \vc{d}_{N,M} : \Grad \vc{v}
\end{equation}
\[
  - \int_\Omega (\vc{f}(\bd_{N,M}) - \Delta \bd_{N,M} )  \otimes \bd_{N,M} : \Grad \vc{v},
\]
\begin{equation}\nonumber
\int_\Omega \vc{u}_{N,M}(0, \cdot) \cdot \vc{v} = \int_\Omega
\vc{u}_0 \cdot \vc{v},
\end{equation}
for any $\vc{v} \in X_N$,
where $r\in(3,10/3)$.
The extra term $ \frac{1}{M} |\Grad\vc{u}_{N,M}|^{r-2}\Grad\vc{u}_{N,M}$ guarantees
sufficient regularity for the velocity field needed in the  director equation.
Our strategy is to pass to the limit first for $N\to\infty$ and then
for $M\to\infty$.

The functions $\vc{d}_{N,M}$ are determined in terms of
$\vc{u}_{N,M}$ as the unique solution of the parabolic system
\begin{equation} \label{approx2}
\partial_t \vc{d}_{N,M}  + \vc{u}_{N,M} \cdot \Grad \vc{d}_{N,M}
- \vc{d}_{N,M}\cdot \Grad \vc{u}_{N,M}  = \Delta \vc{d}_{N,M}- \vc{f}(\vc{d}_{N,M}) ,
\end{equation}
supplemented with
\begin{equation} \label{approx3}
\Grad (d_{N,M})_i \cdot\vc{n}|_{\partial \Omega} = 0,\ i = 1,2,3,
\end{equation}
\begin{equation} \label{approx4}
\vc{d}_{N,M}(0, \cdot) = \vc{d}_{0,M},
\end{equation}
where $\vc{d}_{0,M}$ is a suitable smooth approximation of
$\vc{d}_0$.

Next, given $\vc{u}_{N,M}$, $\vc{d}_{N,M}$, the temperature
$\teta_{N,M}$ is determined as the unique solution to the
heat equation (cf.~Ladyzhenskaya et al.~\cite[Chapter V, Theorem 8.1]{lsu}):
\begin{equation} \label{approx5}
\partial_t \teta_{N,M} + \dive (\teta_{N,M} \vc{u}_{N,M} ) + \dive
\vc{q}^d_{N,M}
\end{equation}
\[=  \tn{S}_{N,M}: \Grad \vc{u}_{N,M} + \frac1M
|\Grad\vc{u}_{N,M}|^r
+\left|\Delta \bd_{N,M}-\vc{f}(\db_{N,M})\right|^2,
\]
\begin{equation} \label{approx6}
\vc{q}^d_{N,M} \cdot \vc{n}|_{\partial \Omega} = 0,
\end{equation}
\begin{equation}\label{approx7}
\teta_{N,M} (0, \cdot) = \teta_{0,M},
\end{equation}
where $\tn{S}_{N,M}=\frac{\mu(\teta_{N,M})}2 \left( \Grad \ub_{N,M} + \Grad^t \ub_{N,M} \right)$, and
\begin{equation}\nonumber
\vc{q}^d_{N,M} = - k( \teta_{N,M} ) \Grad \teta_{N,M} -
h(\teta_{N,M}) \vc{d}_{N,M} (\vc{d}_{N,M} \cdot
\Grad \teta_{N,M} ).
\end{equation}
Actually, relation \eqref{approx5} is just an explicit reformulation
of~\eqref{eqsbis}.

Finally, the pressure $p_{N,M}$ is found as before as
the (unique) solution to a system of integral identities:
\begin{equation} \label{pressap}
\int_\Omega p_{N,M} \Delta \varphi
\end{equation}
\[= \int_{\Omega} \big(\tn{S}_{N,M}
-  \Grad \vc{d}_{N,M} \odot \Grad \vc{d}_{N,M} -(\Delta \bd_{N,M}-\vc{f}(\bd_{N,M}))\otimes \bd_{N,M})-
 [\vc{u}_{N,M}]_M\otimes \vc{u}_{N,M} \big) : \Grad^2 \varphi
\]
\begin{equation}\nonumber
 \mbox{}+ \frac1M \io |\Grad\vc{u}_{N,M}|^{r-2}\Grad\vc{u}_{N,M}:\Grad^2\varphi,
\end{equation}
satisfied
for any test function $\varphi \in C^\infty(\overline{\Omega})$,
$\Grad \varphi \cdot \vc{n}|_{\partial \Omega} = 0$.

Regularizing the convective terms in (\ref{approx1})
is in the spirit of Leray's original approach
\cite{leray} to the Navier-Stokes system.
As a result, we recover the internal energy {\it equality}\
at the level of the limit $N \to \infty$. This fact, in turn, enables us to
replace the internal energy equation (\ref{approx5}) by the
total energy balance before performing the limit $M \to \infty$.
For fixed $M,N$, problem (\ref{approx1}--\ref{pressap}) can be solved
by means of a simple fixed point argument, exactly as in \cite[Chapter 3]{FN}.
Note that all the {\it a priori}\
bounds derived formally in Section \ref{a} apply to our
approximate problem. Thus, given $\vc{u} \in C([0,T]; X_N)$, we can find
$\vc{d} = \vc{d}[\vc{u}]$ solving (\ref{approx2}--\ref{approx4}),
and then $\teta = \teta [\vc{u} , \vc{d} ]$ and the pressure $p$
satisfying (\ref{approx5}--\ref{pressap}). Plugging these functions $\vc{d}$,
$\teta$ in (\ref{approx1}), the corresponding solution ${\cal T}[\vc{u}]$
then defines a mapping $\vc{u} \mapsto {\cal T}[\vc{u}]$.
By the {\it a priori}\ bounds obtained in Section
\ref{a}, we can easily show that ${\cal T}$ possesses a fixed point
by means of the classical Schauder's argument, at least on a possibly short time interval.
However, using once more the {\it a priori} estimates we easily conclude that
the approximate solutions can be extended to any
fixed time interval $[0,T]$, see \cite[Chapter 6]{FN} for details.


\subsection{Passage to the limit as $N\to\infty$}

Having constructed the approximate solutions
$\vc{u}_{N,M}$, $\vc{d}_{N,M}$, $\teta_{N,M}$, and $p_{N,M}$, we
let $N \to \infty$. To take the limit, we need to modify a bit the
formal estimates obtained in Section \ref{a} taking care
of the regularizing terms added in \eqref{approx1} and \eqref{approx5}.
Indeed, from the energy estimate we now additionally obtain
\begin{equation}\label{unuova}
  M^{-1}\|\Grad\vc{u}_{N,M}\|^r_{L^r((0,T)\times\Omega; \RR^{3\times3})}\leq C\,,
\end{equation}
whence we infer that $|\Grad\vc{u}_{N,M}|^{r-2}\Grad\vc{u}_{N,M}$
is uniformly bounded in $L^{\frac{r}{r-1}}((0,T)\times\Omega)$ for fixed $M$.
Moreover, in place of \eqref{apr7} we deduce from \eqref{pressap} the estimate
\begin{equation}\label{pnm}
 \|p_{N,M}\|_{L^{r/r-1}((0,T)\times\Omega)}\leq C(M)\,,
\end{equation}
where we observe that
\begin{equation}\label{rrmu}
  \frac{r}{r-1}\in\Big(\frac{10}7,\frac32\Big),~
    \text{since }r\in\Big(3,\frac{10}3\Big).
\end{equation}

Note that, at least at the level of approximate solutions, relation \eqref{eqteta} holds true as an equality. Hence, taking
$H(\teta)=(1+\teta)^\eta$, with $\eta\in (0,1)$, in \eqref{eqteta},
we get
$$
\|\dt\teta_{N,M}^\nu\|_{(C^0([0,T];W^{1, s}(\Omega)))^*}\leq C \| \dt\teta_{N,M}^\nu\|_{L^1((0,T)\times\Omega)}\leq C,
$$
where $C$ is a positive constant independent of $N$ and $M$, with $s\in (3,+\infty)$, $\nu\in (0,1/2)$.
This leads to the convergence relations:
\begin{align}
\label{cuN}
&\ub_{N,M}\to\ub_M \ \mbox{ weakly-(*) in } L^\infty(0,T;L^2(\Omega;\RR^3))
\cap L^2(0,T;W^{1,2}(\Omega;\RR^3))\,,\\
\label{cuNnuova}
&\Grad\ub_{N,M}\to\Grad\ub_M \ \mbox{ weakly in } L^{r}(0,T;L^r(\Omega;\RR^3))\,,\\
\label{cutN}
& \dt\ub_{N,M}\to\dt\ub_M \ \mbox{ weakly in }
 L^2(0,T;(W^{1,2}(\Omega;\RR^3))^*)
  + L^{\frac{r}{r-1}}(0,T; W^{-1, r/r-1}(\Omega;\RR^3)) \,,\\
\label{cpN}
&p_{N,M}\to p_M \ \mbox{ weakly in } L^{r/r-1}((0,T)\times \Omega)\,,\\
\label{ctetaN}
&\teta_{N,M}^\nu\to\teta_M^\nu\ \mbox{ weakly-(*) in } L^2(0,T;W^{1,2}(\Omega))\cap L^\infty(0,T; L^{1/\nu}(\Omega))\,,\\
\label{ctetatN}
&\dt\teta_{N,M}^\nu\to\dt\teta_M^\nu\ \mbox{ weakly-(*) in } (C_0(0,T;W^{1,s}(\Omega)))^*\,,\\
\label{clambdatetaN}
&\log\teta_{N,M}\to \log\teta_M \mbox{ weakly in } L^2(0,T;W^{1,2}(\Omega))\,,\\
\label{cdN}
&\bd_{N,M}\to \bd_M \ \mbox{ weakly-(*) in } L^\infty(0,T;W^{1,2}(\Omega;\RR^3))
\cap L^2(0,T;W^{2,2}(\Omega;\RR^3))\,,\\
\label{cddtprovv}
& \dt \bd_{N,M}\to \dt \bd_M \ \mbox{ weakly in } L^{5/3}(0,T;L^{5/3}(\Omega;\RR^3))\,.
\end{align}
for any $\nu\in (0,1/2)$, $s>3$, where \eqref{cpN} follows from \eqref{pnm}.
Note that the $M$-projection is kept in the convective term in the limit $N \to \infty$.

Applying the Aubin-Lions compactness
lemma (cf.~\cite{simon}), we deduce that,
\begin{align}
&\teta_{N,M}\to \teta_M \ \mbox{ strongly in } L^p((0,T)\times\Omega)
\end{align}
for any $p\in [1,5/3)$.
Moreover, using \eqref{cdN}, \eqref{cddtprovv}, a simple interpolation argument,
and the Aubin-Lions lemma, we obtain that
\begin{equation}\label{custg}
\Grad\bd_{N,M}\to \Grad \bd_M \ \mbox{ strongly
in }L^\eta((0,T)\times\Omega; \RR^{3\times3})\ \mbox{ for }\eta\in [1,10/3).
\end{equation}
Next, using \eqref{cuN}, \eqref{cuNnuova}, standard interpolation
and embedding properties of Sobolev spaces, and the Aubin-Lions lemma,
we arrive at
\begin{equation}\label{custrN}
  \ub_{N,M}\to\ub_M \ \mbox{ strongly in } L^s((0,T)\times\Omega;\RR^3)\,,\\
\end{equation}
for some $s>5$. Combining this with \eqref{custg}, we finally obtain
\begin{equation}\label{convls}
 \vc{u}_{N,M}\cdot\Grad\bd_{N,M}\to \vc{u}_M\cdot\Grad \bd_M\
 \mbox{ strongly in } L^q((0,T)\times\Omega; \RR^3)
\end{equation}
for some $q > 2$.
Moreover, from \eqref{cuNnuova} and \eqref{cdN}, we have
\[
\bd_{N,M}\cdot \Grad \ub_{N,M} \to \bd_M\cdot\Grad\ub_M \  \mbox{ weakly in } L^p((0,T)\times \Omega; \RR^3)
\]
for some $p>2$, whence
\begin{equation}\label{cdtNbis}
 \dt \bd_{N,M}\to \dt \bd_M \ \mbox{ weakly in } L^2(0,T;L^{2}(\Omega;\RR^3))\,.
\end{equation}
Finally, we have that
\[
|\Grad\ub_{N,M}|^{r-2}\Grad \ub_{N,M}\to \overline{|\Grad\ub_{M}|^{r-2}\Grad \ub_{M}}\
\mbox{weakly in } L^{r/r-1}((0,T)\times \Omega; \RR^{3\times3}).
\]

We conclude that the limit quantities $\vc{u}_M$,
$\vc{d}_M$, $\teta_M$, and $p_M$ solve the problem
\begin{equation}\label{aapprox2}
\int_{\Omega} \vc{u}_M(t, \cdot) \cdot \Grad \varphi = 0 \ \mbox{for a.a.} \ t \in (0,T)
\end{equation}
for any test function $\varphi \in C^{\infty} (\overline{\Omega})$;
\begin{equation} \label{aapprox1}
  \int_0^T \int_{\Omega} \Big( \vc{u}_M \cdot \partial_t \varphi +
   [\vc{u}_{M}]_M\otimes\vc{u}_M : \Grad \varphi \Big) +p_M \dive \varphi= \int_0^T
   \int_{\Omega} (\tn{S}_M+\sigma^{nd}_M) : \Grad \varphi
\end{equation}
\begin{equation} \no
\mbox{}
  - \int_{\Omega} \vc{u}_0
\cdot \varphi (0, \cdot)+ \frac1M \io \overline{| \Grad\vc{u}_{M} |^{r-2} \Grad\vc{u}_{M}}
   : \Grad \varphi,
\end{equation}
for any $\varphi \in C^\infty_0 ([0,T) \times \overline{\Omega}; \RR^3)$,
$\varphi \cdot \vc{n}|_{\partial \Omega} = 0$,
where
\begin{equation}\label{defTM}
  \sigma^{nd}_M = -\left( \Grad \bd_M \odot \Grad \bd_M\right) - \left(\Delta\bd_M-\vc{f}(\bd_M)\right)\otimes \bd_M\,;
\end{equation}
and
\begin{equation}\label{defSM}
\tn{S}_M= \mu( \teta_M) \left( \frac{\Grad \ub_M + \Grad^t \ub_M}2 \right).
\end{equation}

Letting $N \to \infty$ in the equation for $\vc{d}_{N,M}$
we get
\begin{equation} \label{aapprox3}
\partial_t \vc{d}_{M}
+ \vc{u}_{M} \cdot \Grad \vc{d}_{M}-\bd_M\cdot\Grad\ub_M
 = \Delta \vc{d}_{M}-\vc{f}(\vc{d}_{M}),\quad\mbox{a.e. in }(0,T)\times\Omega,
\end{equation}
supplemented with
\begin{equation} \label{aapprox4}
\Grad (d_{M})_i \cdot\vc{n}|_{\partial \Omega} = 0,\ i = 1,2,3,
\end{equation}
\begin{equation} \label{aapprox5}
\vc{d}_{M}(0, \cdot) = \vc{d}_{0,M}.
\end{equation}

The passage to the limit in \eqref{approx5} is more delicate. Actually,
the weak lower semi-continuity of convex functionals on the right-hand side
gives rise to
\begin{equation} \label{aapprox6}
 \partial_t \teta_{M}
 + \dive (\teta_{M} \vc{u}_M )
 + \dive \vc{q}^d_{M}
 \geq
  \frac1M | \Grad\vc{u}_{M} |^r
 + \tn{S}_{M}: \Grad \vc{u}_M
 + \left| \Delta \bd_M - \vc{f}(\bd_M) \right|^2
\end{equation}
satisfied in the sense of distributions,
with
\begin{equation} \label{aapprox7}
  \vc{q}_{M}^d \cdot \vc{n}|_{\partial \Omega} = 0,
\end{equation}
\begin{equation}\label{aapprox8}
\teta_{M} (0, \cdot) = \teta_{0,M},
\end{equation}
where
\begin{equation}\nonumber
\vc{q}^d_{M} = - k( \teta_{M} ) \Grad \teta_{M} -
h(\teta_{M}) \vc{d}_{M} (\vc{d}_{M} \cdot
\Grad \teta_{M} ).
\end{equation}

Next, we claim that the total energy is conserved, namely
\begin{equation}\label{consen}
\partial_t \int_{\Omega} \Big( \frac{1}{2}
|\vc{u}_M|^2 + \teta_M +\frac12|\Grad\bd_M|^2+F(\bd_M)\Big) = 0.
\end{equation}
Indeed, combining \eqref{approx1} with $\vc{v}=\ub_{N,M}$
and (\ref{approx5}--\ref{approx6}), we obtain
\[
  \partial_t\io \left(\frac12|\ub_{N,M}|^2+\teta_{N,M}\right)
\]
\[
  = \io \Big( (\Grad\bd_{N,M}\odot\Grad\db_{N,M}) \cdot \Grad\ub_{N,M}
   +\left((\Delta\bd_{N,M}-\vc{f}(\bd_{N,M}))\otimes \bd_{N,M}\right):\Grad\ub_{N,M}\]
\[
  +|\Delta \bd_{N,M}-\vc{f}(\db_{N,M})|^2\Big),
\]
whence, by virtue of \eqref{approx2} and after a straightforward manipulation, we get
\[
\partial_t \int_{\Omega} \Big( \frac{1}{2}
|\vc{u}_{N,M}|^2 + \teta_{N,M} +\frac12|\Grad\bd_{N,M}|^2+F(\bd_{N,M})\Big) = 0,
\]
yielding, by passing to the limit as $N\to\infty$, the desired conclusion \eqref{consen}.

Now, we want to show that \eqref{aapprox6} is actually an equality.
Taking $\vb =\ub_{N,M}$ in \eqref{approx1} we get
\begin{align}\label{compaMN}
 \|\ub_{N,M}(t)\|^2_{L^2(\Omega)}&
  + \frac12 \int_0^T\io\mu(\teta_{M,N})|\Grad \vc{u}_{N,M} + \Grad^t
  \vc{u}_{N,M}|^2
+\frac2M\int_0^T\io |\Grad\vc{u}_{N,M}|^r\\
\no
&=\|\ub_0\|^2_{L^2(\Omega)}+2\int_0^T\io\ \sigma^{nd}_{N,M}
: \Grad \vc{u}_{N,M}\,.
\end{align}

Next,
thanks to  (\ref{cuN}--\ref{cutN}),
we can
take $\bu_M$ as a test function in \eqref{aapprox1}
to obtain
\begin{align}\label{compaM}
 \|\ub_M(t)\|^2_{L^2(\Omega)}&
  +\frac12  \int_0^T\io\mu(\teta_M)|\Grad \vc{u}_{M} + \Grad^t
\vc{u}_{M}|^2\\
 \nonumber
  & \mbox{} ~~~~~~~~
  + \frac2M\int_0^T\io \overline{|\Grad\vc{u}_{M}|^{r-2}\Grad\bu_M}: \Grad \bu_M \\
 \nonumber
  & = \|\ub_0\|^2_{L^2(\Omega)}+2\int_0^T\io \sigma^{nd}_M
   : \Grad \vc{u}_{M}\,.
\end{align}
Now, multiplying \eqref{approx2} by $\Delta\bd_{N, M} -\vc{f}(\bd_{N,M})$, we obtain
\begin{equation}\label{compabisNM}
\|\Grad\bd_{N,M}(t)\|^2_{L^2(\Omega)}+\io F(\bd_{N,M})(t)+2\int_0^T\io|\Delta\bd_{N,M}-\vc{f}(\bd_{N,M})|^2
\end{equation}
\[
=\|\Grad\bd_0\|^2_{L^2(\Omega)}+2\io F(\bd_0)+2\int_0^T(\ub_{N,M}\cdot\Grad\bd_{N,M}-\bd_{N,M}\cdot\Grad\ub_{N,M},\Delta\bd_{N,M}-\vc{f}(\bd_{N,M})).
\]
Analogously, multiplying \eqref{aapprox3} by $\Delta \bd_{M} -\vc{f}(\bd_M)$ we get
\begin{equation}\label{compabisM}
\|\Grad\bd_{M}(t)\|^2_{L^2(\Omega)}+\io F(\bd_{M})(t)+2\int_0^T\io|\Delta\bd_{M}-\vc{f}(\bd_{M})|^2
\end{equation}
\[
=\|\Grad\bd_0\|^2_{L^2(\Omega)}+2\io F(\bd_0)+2\int_0^T(\ub_{M}\cdot\Grad\bd_{M}-\bd_{M}\cdot\Grad\ub_{M},\Delta\bd_{M}-\vc{f}(\bd_{M})).
\]
Taking the sum of \eqref{compaMN} and \eqref{compabisNM} and \eqref{compaM} with \eqref{compabisM}, and,
finally, passing to the limit as $N\to\infty$,
we obtain
\[
\int_0^T\io |\Grad \ub_{N,M}|^r\to \int_0^T\io \overline{|\Grad\ub_{M}|^{r-2}\Grad\ub_M}:\Grad\ub_{M},
\]
\[
\int_0^T\io|\Delta\bd_{N,M}-\vc{f}(\bd_{N,M})|^2\to \int_0^T\io|\Delta\bd_{M}-\vc{f}(\bd_{M})|^2,
\]
entailing, by means of standard Minty's trick and monotonicity argument,
\[
\Grad\ub_{N,M}\to\Grad\ub_M\ \mbox{strongly in }L^{r}((0,T)\times\Omega;\RR^{3\times 3}),
\]
\[
\Delta \bd_{N,M}\to \Delta \bd_M\ \mbox{strongly in }L^{2}((0,T)\times\Omega;\RR^3).
\]

Consequently, the inequality \eqref{aapprox6} may be replaced by the equality
\begin{equation} \label{aapprox6+}
 \partial_t \teta_{M} + \dive (\teta_{M} \vc{u}_M ) + \dive
\vc{q}^d_{M}
=
 \frac1M | \Grad\vc{u}_{M} |^r
+ \tn{S}_{M}: \Grad \vc{u}_M +\left|\Delta \bd_{M}-\vc{f}(\bd_M)\right|^2.
\end{equation}

Taking $\ub_M\varphi$, with $\varphi\in {\cal D}((0,T)\times\Omega)$,
as a test function in \eqref{aapprox1},
testing \eqref{aapprox3} by $\frac{D \bd_M}{Dt}\varphi$,
adding both relations to \eqref{aapprox6+} multiplies by $\varphi$,
and using \eqref{edot}, we get an $M$-analogue of \eqref{weak4}, namely:
\begin{equation}\label{totaleM}
\partial_t \left( \frac{1}{2}|\bu_M|^2 + e_M \right)
+ \dive \left( \frac{1}{2}|\bu_M|^2[\ub_M]_M + e_M \ub_M \right)
\end{equation}
\[
+ \dive \Big( p_M \bu_M + \vc{q}_M - \tn{S}_M \bu_M - \sigma^{nd}_M \bu_M \Big)
= \dive \Big( \Grad \bd_M\cdot \left( \Delta \bd_M - \vc{f}(\bd_M) \right) \Big),
\]
with the internal energy
\[
e_M = \frac{1}{2} |\Grad \bd_M |^2 + F(\bd_M) + \theta_M
\]
and the flux
\[
{\bf q}_M  = -k( \theta_M )\Grad\theta_M -h(\theta_M )(\db_M\cdot \Grad\theta_M )\db_M
  - \lambda \Grad \bd_M \cdot \Grad \bu_M \cdot \bd_M.
  \]

Finally,  we can multiply \eqref{approx5} by $H'(\teta_M)\varphi$, obtaining
\begin{equation}\label{weak5M}
  \int_0^T\io H(\teta_M) \partial_t\varphi +\int_0^T\io \left(H(\teta_M) \ub_M
   + H'(\teta_M)\vc{q}_M^d\right)\cdot\Grad\varphi
\end{equation}
\begin{align*}
  & \leq -\int_0^T\io \bigg(H'(\teta_M)\Big(\tn{S}_M: \Grad \bu_M
   + \frac{1}M | \Grad \bu_M |^r \\
  & \mbox{}~~~~~
   + | \Delta \bd_M - \vc{f}(\bd_M) |^2 \Big)
   + H''(\theta_M) \vc{q}_M^d \cdot \Grad \theta_M\bigg)\varphi
\end{align*}
\[
  -\io H(\teta_{0, M})\varphi(0,\cdot),
\]
for any $\varphi \in C^{\infty}_0 ([0,T) \times \overline{\Omega})$,
$\varphi\geq 0$, and any smooth, non-decreasing and concave function
$H$. To be precise, we have however to remark that, at this level,
we do not have sufficient regularity in \eqref{approx5} to use
$H'(\teta_M)\varphi$ directly as a test function. Nevertheless, the procedure
could be justified by a standard regularization argument and then taking
the (supremum) limit. This is also the reason why we get the $\le$ sign,
rather than the equality, in \eqref{weak5M}. This concludes the
passage to the limit for $N\to\infty$.


\subsection{Passage to the limit as $M\to\infty$}

Our final goal is to let $M\to\infty$
in (\ref{aapprox2}--\ref{aapprox5}), \eqref{totaleM}, and \eqref{weak5M}.
We notice that the limits in (\ref{cuN}), (\ref{ctetaN}--\ref{custg})
still hold  when letting  $M\to \infty$.
On the other hand, we now have
\begin{align}
\label{cdtuM}
& \dt\ub_M\to\dt\ub\ \mbox{ weakly in } L^{\frac{r}{r-1}}(0,T;W^{-1,\frac{r}{r-1}}(\Omega;\RR^3))\,,\\
\label{cpM}
&p_M\to p\ \mbox{ weakly in } L^{\frac{r}{r-1}}((0,T)\times\Omega)\,,\\
\label{cdtM}
&\dt\bd_{M}\to \dt \bd \ \mbox{ weakly in } L^2(0,T;L^{3/2}(\Omega;\RR^3))\,,
\end{align}
and, obviously,
\begin{equation}
\label{cduMM}
M^{-1/(r - 1)}\Grad\ub_M\to 0 \mbox{ strongly in } L^{r - 1}((0,T)\times \Omega).
\end{equation}

The above relations are sufficient to
pass to the limit $M \to \infty$ in (\ref{aapprox2}--\ref{aapprox5}) to recover
(\ref{weak1}--\ref{weak3}). In addition, by \eqref{rrmu} and the previous estimates,
we get
\begin{align}\nonumber
&\left\{ \left(\frac{|\ub_M|^2}{2}+p_M\right)\ub_M \right\}_{M > 0}
\ \mbox{bounded in } L^{\iota}((0,T)\times\Omega; \RR^3)\ \mbox{for some}\ \iota>1\,,\\
\nonumber
& \left\{ \teta_M \ub_M \right\}_{M > 0} \ \mbox{bounded in } L^{q}((0,T)\times\Omega;\RR^3))\ \mbox{for any}\ q\in [1, 10/9)\,,\\
\nonumber
& \left\{\sigma^{nd}_M\ub_M \right\}_{M > 0}\mbox{bounded in } L^{\iota}((0,T)\times\Omega; \RR^3)\ \mbox{for some}\ \iota>1 \,,\\
\nonumber
&{\bf q}_M\mbox{ bounded in }\ \mbox{bounded in } L^{\iota}((0,T)\times\Omega; \RR^3)\ \mbox{for some}\ \iota>1\,.
\end{align}
Notice that we used here in an essential way the fact that $r/(r-1)>10/7$.

As a consequence, we can pass to the limit in \eqref{totaleM} and to the $\limsup$ in \eqref{weak5M}
(thanks also to the positivity and convexity
of the terms on the last line of \eqref{weak5M}) to deduce the desired conclusions \eqref{weak4} and
\eqref{weak5}. This completes the proof of Theorem
\ref{theo1}.

To conclude, we remark that the above estimates {\em are not} sufficient
for passing to the limit in \eqref{aapprox6+}
with respect to $M\to\infty$, due to the lack
of strong convergences of the terms appearing on the
\rhs.
%
%


\end{document}